\documentclass[12pt]{article}
\usepackage[final]{epsfig}
\usepackage{graphics}
\usepackage{amsmath}
\usepackage{amsfonts}
\usepackage{latexsym}
\usepackage{amssymb}
\usepackage{graphicx}
\usepackage{url}
\usepackage{epstopdf}
\usepackage{hyperref}
\usepackage{color}
\usepackage{marginnote}
\usepackage{a4wide}

\newtheorem{lemma}{Lemma}[section]

\newtheorem{theorem}{Theorem}

\newcommand{\g}{{\gamma}}

\newcommand{\proofend}{$\Box$\bigskip}

\newcommand{\R}{{\mathbb R}}

\newcommand{\RP}{{\mathbb {RP}}}

\def\proof{\paragraph{Proof.}}

\begin{document}

\title{On cusps of caustics by reflection in  two dimensional projective Finsler metrics}

\author{
Serge Tabachnikov
 \footnote{
Department of Mathematics,
Pennsylvania State University,
University Park, PA 16802,
USA;
tabachni@math.psu.edu}
}

\date{\today}

\maketitle

\section{Motivation and previous results} \label{sect:mot}

In the posthumously published ``Lectures on Dynamics", Jacobi claimed that the conjugate locus of a non-umbilic point of a triaxial ellipsoid has exactly four cusps. This is known as the Last Geometric Statement of Jacobi. The conjugate locus of a point is the envelope of the geodesics that emanate from this point. These geodesics have the second, third, etc., envelopes; they are also called the first, second, etc., caustics. 

The Last Geometric Statement of Jacobi was proved relatively recently  \cite{IK}.  Conjecturally, each next caustic also has exactly four cusps, see \cite{Si}.  One also has a theorem, attributed to C. Carath\'eodory by W. Blaschke in his differential geometry textbook: 
 {\it The conjugate locus of a generic point on a convex surface has at least four cusps}. See \cite{Wa} for a recent proof.

One may consider a billiard version of this problem: instead of a closed surface, take a billiard table in the Euclidean plane bounded by an oval (smooth strictly convex closed curve), and instead of the pencils of geodesics, consider the pencil of billiard trajectories starting at a point inside the billiard table.  After $n$ reflections off the boundary, one obtains a 1-parameter family of lines, and their envelope is the $n$th {\it caustic by reflection}. One may use the language of geometrical optics: the point is a source of light and the boundary curve is an ideal mirror. 

\begin{figure}[ht]
\centering
\includegraphics[width=.4\textwidth]{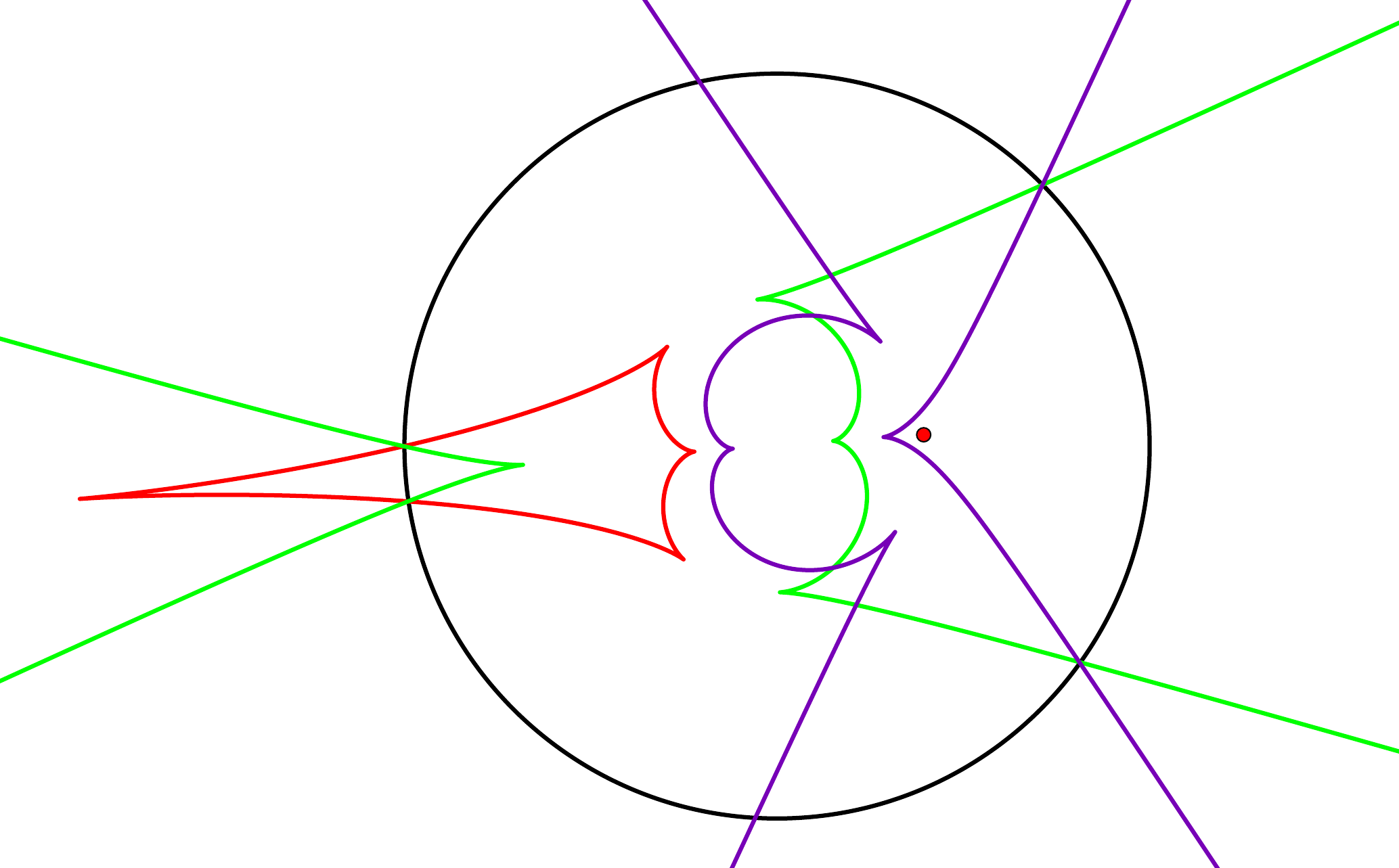}
\caption{First three caustics by reflection in a circle.}
\label{three}
\end{figure}

This billiard problem was studied in two recent papers \cite{BT,BST}. We proved that for every oval, every $n\ge 1$, and a generic source of light the $n$th caustic by reflection has at least four cusps. We provided some evidence toward the conjecture that this number is exactly four for all $n$  if the billiard table is elliptic and proved this conjecture in the case when the boundary curve is a circle, see Figure \ref{three}. The  context for these results is the famous 4-vertex theorem and its numerous variations and generalizations; see, e.g., \cite{Ar}. 

In this note we extend this four cusps result to Finsler billiards in the special case of a projective Finsler metric, a (not necessarily symmetric)  Finsler metric whose geodesics are straight lines. 

\section{Finsler metrics and Finsler billiards} \label{sect:Fin}

From the point of view of geometrical optics, Finsler geometry describes the propagation of light in an inhomogeneous anisotropic medium $M$:  the velocity of light depends on the point and the direction. 

As usual, one has two descriptions of this process, the  Lagrangian and the Hamiltonian ones. From the 
Lagrangian perspective, Finsler metric is defined by a field of {\it indicatrices} $I_x \subset T_x M$, the unit sphere subbundle of the tangent bundle of $M$. These indicatrices are unit level hypersurfaces of the Lagrangian function on $TM$ that defines the metric. 
The indicatrices are smooth and strictly convex hypersurfaces but, in general, not necessarily origin-symmetric.

The dual Hamiltonian description provides a field of {\it figuratrices} $J_x \subset T^*_x M$, the unit cosphere subbundle of the cotangent bundle of $M$. The indicatrices and figuratrices are related by the  Legendre transform $D: I_x \to J_x$ (the polar duality):
$$
I_x \ni v \mapsto w \in J_x\ \ {\rm if}\ \   {\rm Ker}\ w = T_v I_x\ \ {\rm and}\  \ w(v)=1.
$$ 

A Finsler geodesic is a curve that extremizes the Finsler length (or optical path length) between its endpoints. The Finsler geodesic flow  is defined similarly to the Riemannian case: the foot point of a Finsler unit tangent vector moves with the unit speed along the Finsler geodesic that it defines, and the vector remains  unit and tangent to this geodesic. 

We refer to any of the numerous textbooks on Finsler geometry or to the surveys \cite{Al2,Ch}.

Finsler billiard reflection was defined in \cite{GT} similarly to the usual, Riemannian one, by a variational principle.  Let $M$ be a Finsler manifold with boundary $S$, a billiard table,  let $a$ and $b$ be two points in $M$ and $x$ be a boundary point. One says that the Finsler geodesic ray $ax$ reflects to the ray $xb$ if $x$ is a critical point of the Finsler distance function $F(x)={\rm dist} (ax)+{\rm dist}(xb)$. Note that, in general, the reflection is not reversible: it is not necessarily true that $bx$ reflects to $xa$. 

Finsler billiard reflection can be described geometrically. This description is especially nice in dimension 2, which is the subject of this note. We continue to refer to \cite{GT}.

\begin{figure}[ht]
\centering
\includegraphics[width=.6\textwidth]{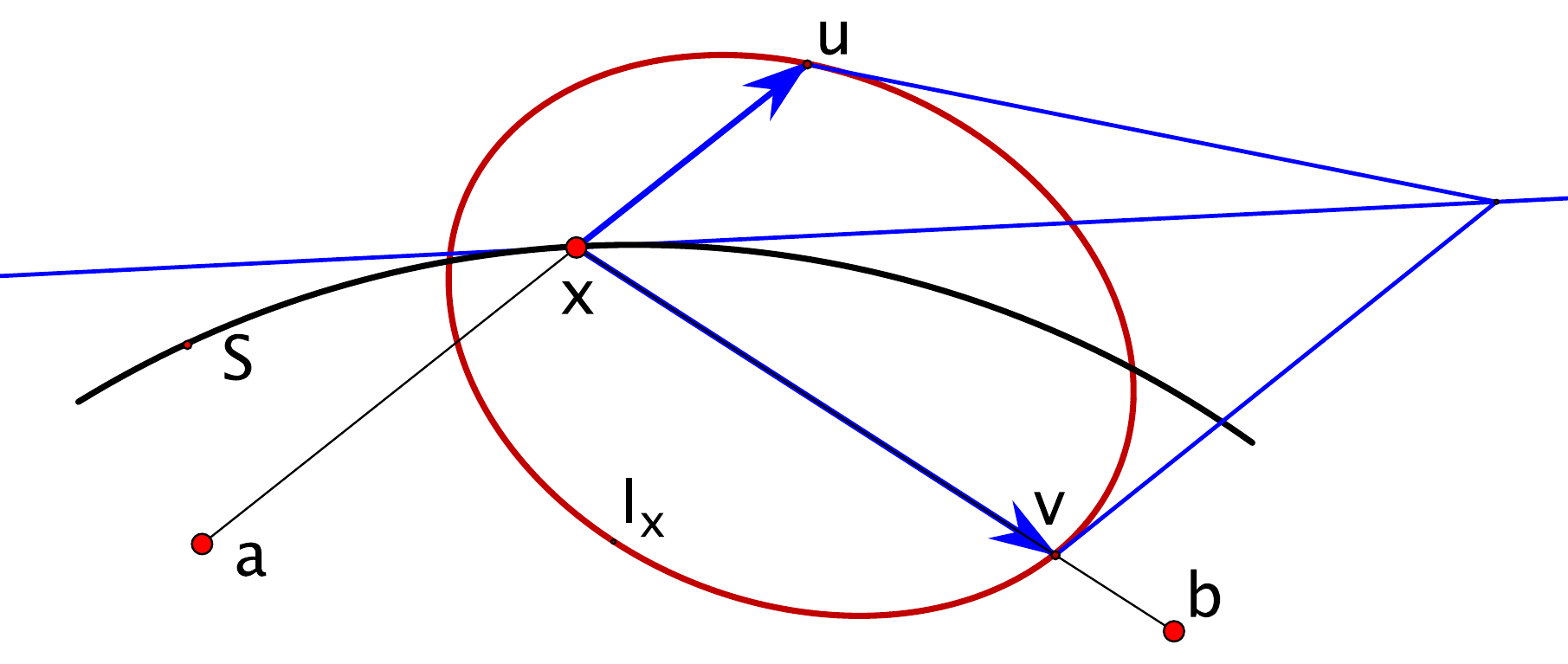}
\caption{Finsler billiard reflection.}
\label{reflection}
\end{figure}

Consider Figure \ref{reflection}.  The red oval is the indicatrix at the reflection point $x \in S$, and $u$ and $v$ are the incoming and outgoing unit velocity vectors. The reflection law states that the tangent lines to the indicatrix at points $u$ and $v$ and the tangent line to the boundary $S$ at point $x$ are concurrent (this includes the case when  the three lines are parallel). 

In the Euclidean case, the indicatrix is a circle, and this reflection law becomes the familiar ``the angle of incidence equals the angle of reflection".

A popular example of a Finsler billiards is a Minkowski billiard. In Minkowski geometry the indicatrices are parallel translation copies of each other and the geodesics  are straight lines. A Minkowski billiard is defined by two ovals, the indicatrix and the billiard table. These billiards were studied in connection with the Viterbo conjecture in symplectic topology and its relation with the Mahler conjecture in convex geometry, see \cite{AKO}.

Concerning Minkowski billiards, also see \cite{Ra}. 

\begin{figure}[ht]
\centering
\includegraphics[width=.36\textwidth]{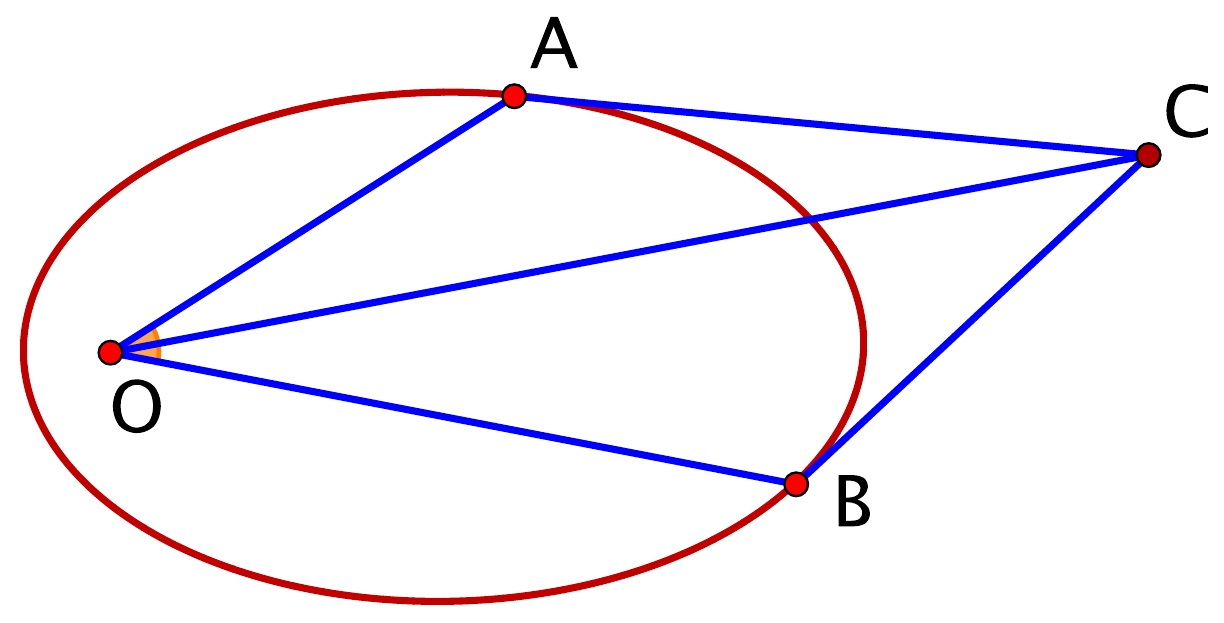}
\caption{The indicatrix is a Kepler ellipse.}
\label{ellipse}
\end{figure}

Consider another example: the indicatrix is a focus-centered (Kepler) ellipse, see Figure \ref{ellipse}. A theorem of elementary geometry states that  $\angle AOC = \angle BOC$, see \cite{AZ}, Theorem 1.4 (this  is known as ``Le second th\'eor\`eme de Poncelet", see the  Wikipedia page ``Th\'eor\`eme de Poncelet"). 

This result implies that the respective Finsler billiard reflection satisfies the same law of equal angles as the usual, Euclidean, one. This applies to another popular billiard model, the magnetic billiards. We follow the  discussion in \cite{Ta}. 

A magnetic field exerts a force on a moving charge that is perpendicular to the direction of motion and is proportional to the speed (Lorentz force).
In particular, the speed of the charge remains constant. 

 A magnetic field in the plane is a given by a differential 2-form $B(x_1,x_2) dx_1 \wedge dx_2$, where the function $B$ is the strength of the magnetic field. Choose a differential 1-form $\alpha$ such that $d\alpha=- B(x_1,x_2) dx_1 \wedge dx_2$. The Lagrangian for the motion of a charge in this magnetic field is  
$$
L(x,v) = \frac{1}{2} |v|^2 + \alpha(x) (v)
$$
(the 1-form $\alpha$ is not unique, but the freedom of its choice does not affect the dynamics).

Following the Maupertuis principle, one replaces the Lagrangian by a homogeneous of degree one Lagrangian 
\begin{equation} \label{eq:Lagr}
L(x,v) =  |v| + \alpha(x) (v)
\end{equation}
whose extremals are the trajectories of the charge moving with the unit speed. In particular, if the magnetic field is constant, 
$$
L(x,v) =  |v| + \frac{1}{2R} \det(v,x),
$$
and the trajectories are counterclockwise oriented circles of (Larmor) radius $R$.

We assume that the magnetic field is sufficiently weak, so that $L(x,v) > 0$ for $v \neq 0$  in the domain under consideration. More precisely, we assume that $|\alpha(x)| < 1$ for all $x$. 

\begin{lemma} \label{lm:ell}
Let $L(x,v)$ be as in (\ref{eq:Lagr}). Then, for every $x$, the indicatrix $L(x,v)=1$ is a focus-centered ellipse.
\end{lemma}

\proof
The equation of a focus-centered axes-aligned ellipse in the $(v_1,v_2)$-plane is
\begin{equation} \label{eq:ellf}
\frac{(v_1+c)^2}{a^2} + \frac{v_2^2}{b^2}=1\ \ {\rm with}\ \ a^2=b^2+c^2.
\end{equation}
We need to show that the equation $L(x,v)=1$ has this form.

Rotating the $v$-plane if needed, we may assume that $\alpha(x)(v)=tv_1$ with $|t|<1$. Then
$$
\sqrt{v_1^2+v_2^2} + tv_1=1,
$$
which is rewritten as
$$
(1-t^2)^2 \left(v_1 + \frac{t}{1-t^2}\right)^2 + (1-t^2) v_2^2 =1.
$$
Therefore, setting
$$
a=\frac{1}{1-t^2},\  b=\frac{1}{\sqrt{1-t^2}},\  c=\frac{t}{1-t^2}
$$
yields the desired equation (\ref{eq:ellf}). 
\proofend

Thus the indicatrices of the magnetic Finsler metric are Kepler ellipses.

\begin{figure}[ht]
\centering
\includegraphics[width=.33\textwidth]{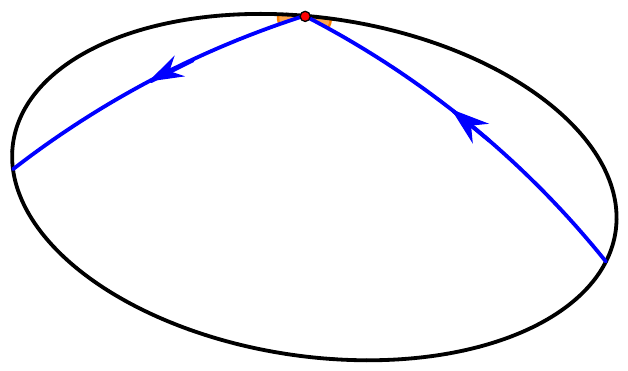}
\caption{Magnetic billiard: the angle of incidence equals the angle of reflection.}
\label{magnetic}
\end{figure}

Magnetic billiards model the the motion of a charge in a  a magnetic field with specular reflection off the boundary of the domain, so that the angle of incidence equals the angle of reflection, see Figure \ref{magnetic}. Due to the  ``second theorem of Poncelet", Figure \ref{ellipse}, 
and Lemma \ref{lm:ell}, magnetic billiards are  a particular case of  Finsler billiards whose metric is defined by the Lagrangian (\ref{eq:Lagr}).

We cannot help mentioning a multidimensional generalization of these results due to Akopyan and Karasev \cite{AK}.

\begin{theorem} \label{thm:AK}
Let $K$ be a smooth convex body in $\R^n$ containing the origin, and let $K_1$ be its convex image under a projective transformation that preserves, with orientation, every line passing through the origin. 
Then the Minkowski billiard reflection law in the space with the norm $K$ is the same as in the space with the norm $K_1$.
\end{theorem}

Such maps are given by the formula
$$
x \mapsto \frac{tx}{1 + \ell(x)},\ t > 0,
$$
where $\ell$ is a linear function. They send  the origin-centered spheres to the focus-centered ellipsoids, in particular, origin-centered circles to Kepler ellipses.

\medskip
Let us finish this section by adressing the symplectic properties of Finsler billiards. This properties do not play a major role in the present note, so we will be brief.
 
 Assume that the space of oriented non-parameterized Finsler geodesics of $M$ is a smooth manifold. The Finsler billiard reflection defines a transformation of this space,  the billiard ball map. This space of geodesics carries a symplectic structure constructed as follows. 
 
 Identify tangent and cotangent vectors via the Legendre transform. The cotangent bundle $T^* M$ has a canonical symplectic structure, and its restriction to the unit cosphere bundle has a 1-dimensional kernel at every point. The integral curves of this field of directions, the characteristics, are identified with the oriented non-parameterized Finsler geodesics of $M$. As a result, the space of characteristics carries a symplectic structure obtained from that in $T^* M$ by restriction to the unit cosphere bundle and factoring out the kernel. This construction is familiar in the Riemannian case, but it extends without change to the Finsler one.
 
 A fundamental feature of Finsler billiards is that the billiard ball map preserves the symplectic structure of the space of oriented non-parameterized geodesics. It is important that this invariant symplectic form does not depend on the shape of the billiard table, it is determined by the ambient Finsler metric only. We refer to \cite{GT} for details. 
 
\section{Projective Finsler metric in two dimensions} \label{sect:two}
In this section, we mostly follow the exposition in \cite{Al}. See \cite{Bu,Po} for more details. 

Hilbert's fourth problem asks to {\it construct and study the geometries in which the straight line segment is the shortest connection between two points}. We interpret this problem (in dimension two,  which is our concern in this note) as asking to describe Finsler metrics in convex subsets of the plane whose geodesics are straight segments. Such metrics are called projective. 

The first examples are provided by Riemannian metrics of constant curvature, the Euclidean, spherical, and hyperbolic ones. The Euclidean case needs no explanation. 

Consider a round sphere and project it from the center to a plane. This central projection takes great circles to straight lines, and it defines a projective metric in the plane that has a constant positive curvature. 

A similar construction works for the hyperbolic plane presented by a hyperboloid of two sheets in Minkowski space: the central  projection takes one sheet of the hyperboloid to the open unit disc, sending  geodesics to straight lines. This yields the projective (Beltrami-Caley-Klein) model of the hyperbolic plane.

According to a Beltrami theorem, a projective Riemannian metric is a metric of constant curvature, so the above examples exhaust the Riemannian  cases.

\begin{figure}[ht]
\centering
\includegraphics[width=.3\textwidth]{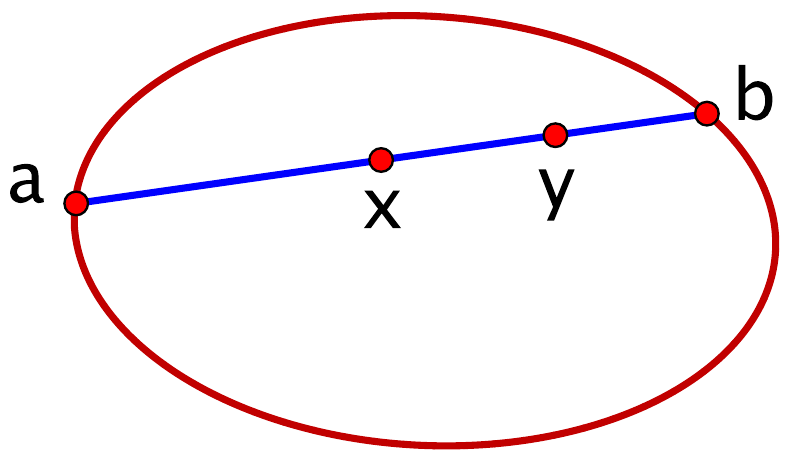}
\caption{Hilbert's and Funk's metrics.}
\label{Funk}
\end{figure}

A Minkowski metric  is an example of a projective Finsler metric. The projective model of the hyperbolic plane generalizes to Hilbert's and to Funk's metrics in a convex domain in the projective plane, see Figure \ref{Funk}, given by the formulas
$$
d_H(x,y)=\frac{1}{2}\ln\left(\frac{|y-a||x-b|}{|y-b||x-a|}\right),\ \ d_F(x,y)=\ln\left(\frac{|x-b|}{|y-b|}\right)
$$
(Hilbert's metric is symmetric, while Funk's metric is not). See \cite{Fa} for a study of Funk billiards.

\begin{figure}[ht]
\centering
\includegraphics[width=.4\textwidth]{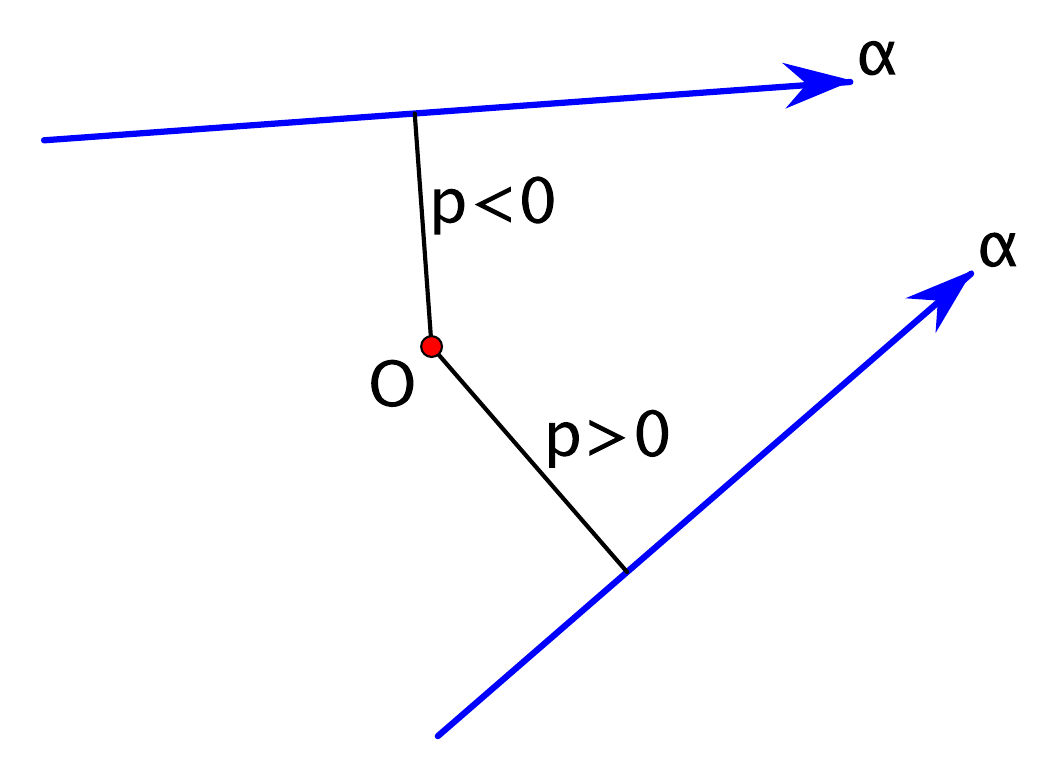}
\caption{Coordinates in the space of oriented lines.}
\label{lines}
\end{figure}

Let us introduce coordinates in the space ${\mathcal L}$ of oriented lines in $\R^2$. Choose an origin $O$. An oriented line is determined by its direction $\alpha \in S^1$ and its signed distance $p \in \R$ from the origin, see Figure \ref{lines}. Thus ${\mathcal L}$ is an infinite cylinder.

The symplectic structure on ${\mathcal L}$, invariant under the billiard ball transformation and described at the end of Section \ref{sect:Fin}, is given by the 2-form $dp \wedge d\alpha$. Up to a factor, this is the unique area form on the space of lines that  is invariant under isometries of the plane. We denote by $dA$ the respective area element.

We briefly describe a construction of symmetric projective Finsler metrics due to H. Buseman. 

Recall the Cauchy-Crofton formula. Let $\g\subset \R^2$ be a piecewise smooth curve. Define a piecewise constant function on ${\mathcal L}$ to be the number of intersection points of a line with the curve. Then
$$
{\rm length} (\g) = \frac{1}{4} \int_{\ell \in \mathcal L} \#(\ell \cap \g)\ dA.
$$

Let $f: {\mathcal L} \to \R$ be a positive smooth function. Replace the area element $dA$ with $f dA$; then an analog of the Cauchy-Crofton formula defines a symmetric projective Finsler metric in the plane. We refer to \cite{Al2} for more information.

\section{The four cusps theorem} \label{sect:four}

We finally turn to the main result of this note. 

Let $U$ be an open plane domain with a projective (not necessarily symmetric) Finsler metric, and let $\g \subset U$ be an oval, the boundary of a Finsler billiard table. Let $O$ be a point inside $\g$ (the source of light). Consider the 1-parameter family of billiard trajectories, starting at $O$ and undergoing $n$ Finsler billiard reflections.

\begin{theorem} \label{thm:main}
For every $n \ge 1$, the envelope of this 1-parameter family of lines in $\RP^2$ (the $n$th caustic by reflection) has at least four cusps.  
\end{theorem}

We need to comment on this formulation. Taken literally, it assumes that the $n$th caustic by reflection is a piecewise smooth curve whose smooth arcs connect distinct generic (semi-cubic) cusps, that is, this caustic is sufficiently generic. Of course, the caustic may degenerate, even to a point: for example, this is the case if point $O$ is a focus of an ellipse in the Euclidean plane.\footnote{One can construct an analog of ellipse in Finsler geometry as the locus of points whose sum of distances to two fixed points is fixed. Such a curve shares the optical property of ellipse.} 

To include possibly degenerate caustics, one can reformulate the statement of the theorem as follows: {\it there exist at least four distinct oriented lines through point $O$ that, after $n$ Finsler billiard reflections, pass through singular points of  the $n$th caustic by reflection.} 

This is similar to the classic 4-vertex theorem: one common formulation is that the evolute of an oval (the envelope of its normals) has at least four cusps, but a more precise statement is that the curvature of this oval has at least four critical points. We prefer the former formulation as being more graphical.

\paragraph{Proof of Theorem.} The phase space of the billiard ball map is the subset of ${\mathcal L}$ consisting of the lines that intersect the curve $\g$; the boundary of this phase cylinder are the two curves comprising the oriented lines tangent to $\g$. We think of ${\mathcal L}$ as the vertical cylinder in $\R^3$ whose axis passes through the origin.  

The space ${\mathcal L}$ carries a 2-parameter family of curves comprising the  lines passing through fixed points. Using the language of the projective duality, we call these curves ``lines". In the $(\alpha,p)$-coordinates, such a ``line" is the sine curve 
$$
p=a\sin\alpha-b\cos\alpha,
$$
where $(a,b)$ is the respective point. These ``lines" are the intersections of the cylinder ${\mathcal L}$ with the planes through the origin.

Let $C_n \subset {\mathcal L}$ be the curve consisting of the lines that  started at point $O$ and made $n$ Finsler billiard reflections. 
This curve is projectively dual to the $n$th caustic by reflection. The cusps of the $n$th caustic by reflection correspond to the second-order tangencies of the curve $C_n$ with ``lines". These are ``inflections" of $C_n$.

The curve $C_n$ goes around the phase cylinder once. Indeed, this is true for the original pencil of lines through point $O$, and hence for its consecutive images under the billiard ball map. 

Consider the central projection of the phase cylinder to the unit sphere. The ``lines" become great circles, and the ``inflections" of the curve $C_n$
become the spherical inflections of its projection to $S^2$. We need to show that there are at least four such inflections.

We use a theorem of B. Segre that {\it if a simple closed spherical curve intersects every great circle, then it has at least four inflection points}, see \cite{Se}.\footnote{This implies Arnold's ``tennis ball theorem": a simple smooth closed spherical curve that bisects the area has at least four spherical inflections.}  Thus we need to show that $C_n$ intersects every ``line".

If a point $A$ lies inside the billiard table, this asserts the existence of an $n$-bounce Finsler billiard shot from $O$ to $A$. Consider $n$ points $x_1,x_2,\ldots,x_n \in \g$, and let $F(x_1,\ldots,x_n)$ be the Finsler length of the polygonal path $O x_1\ldots x_n A$. This function has a maximum, and due to the triangle inequality, at this maximum one has $x_i \ne x_{i+1}$ for all $i$. Hence the polygonal path $O x_1\ldots x_n A$ is the desired billiard trajectory.

\begin{figure}[ht]
\centering
\includegraphics[width=.55\textwidth]{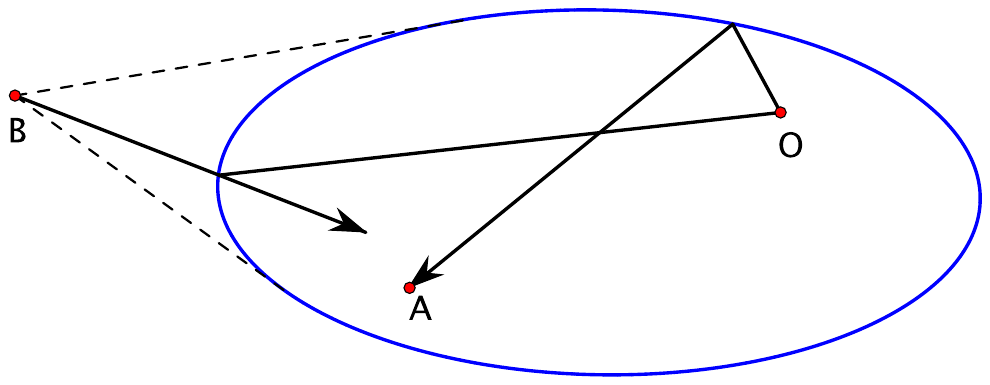}
\caption{The curve $C_n$ intersects every "line".}
\label{intersect}
\end{figure}

If a point $B$ lies outside or on the boundary of the billiard table, the statement holds for a topological reason: the ``line", dual to point $B$, connects the two boundaries of the phase cylinder, and it must intersect the non-contractable curve $C_n$. See Figure \ref{intersect}. 

This concludes the proof (which is a variation of one of the arguments given in \cite{BT}).

Let us remark that there exist at least $n+1$ $n$-bounce billiard shots from $O$ to $A$; this follows from a slight modification of a theorem of M. Farber \cite{Far}. 
\medskip

Let us finish with a problem: {\it Does a 4-cusp result, similar to Theorem \ref{thm:main}, hold for more general Finsler billiards?} 

\begin{figure}[ht]
\centering
\includegraphics[width=.35\textwidth]{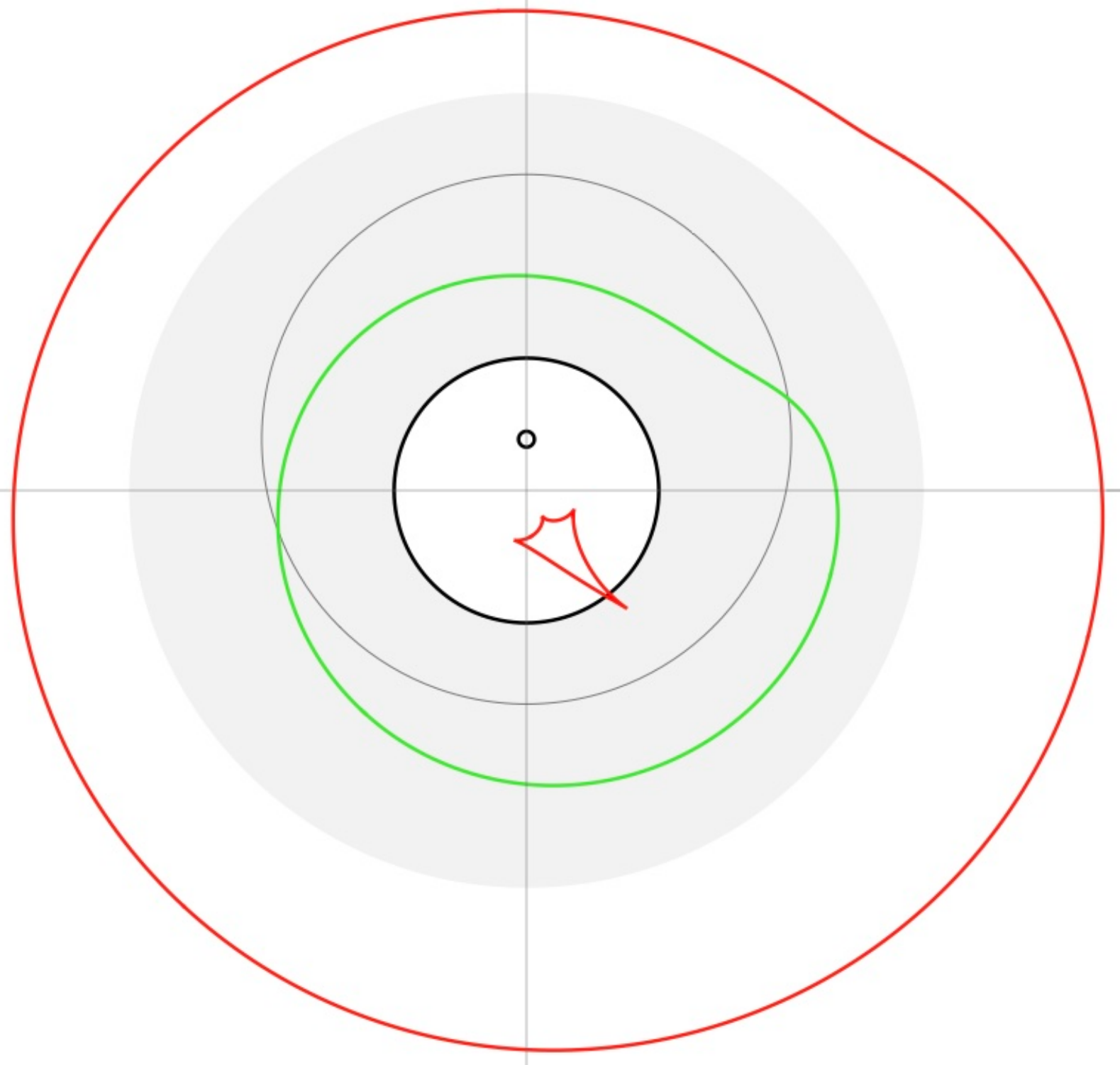} \qquad
\includegraphics[width=.35\textwidth]{circle1}
\caption{The first and second caustics by reflection (courtesy of G. Bor). The inner component of the caustic has four cusps, while the outer one is smooth. The green curve is the curve of centers of the Larmor circles, an analog of the curve $C_n$.}
\label{circles}
\end{figure}

\begin{figure}[ht]
\centering
\includegraphics[width=.35\textwidth]{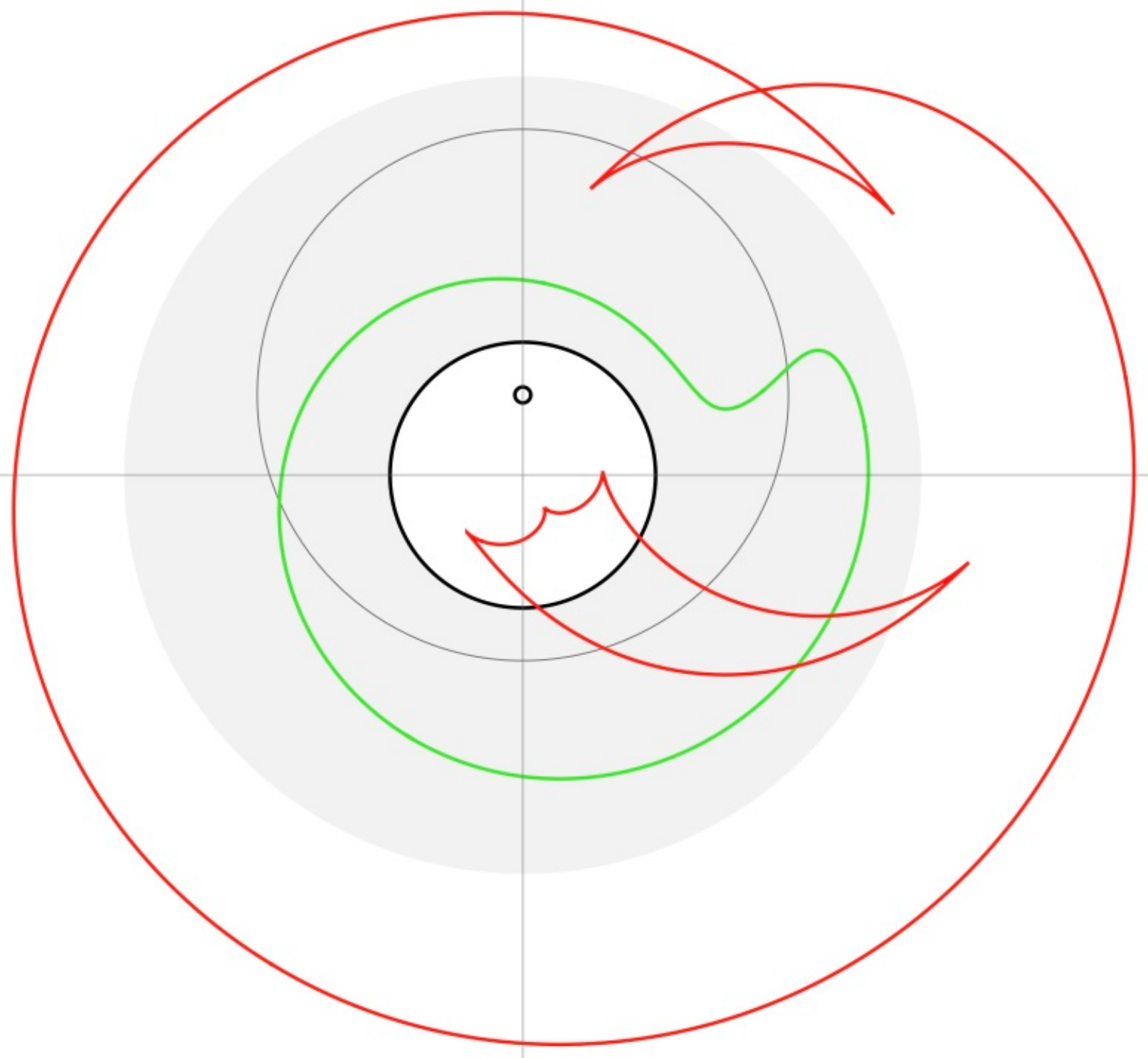} \qquad
\includegraphics[width=.35\textwidth]{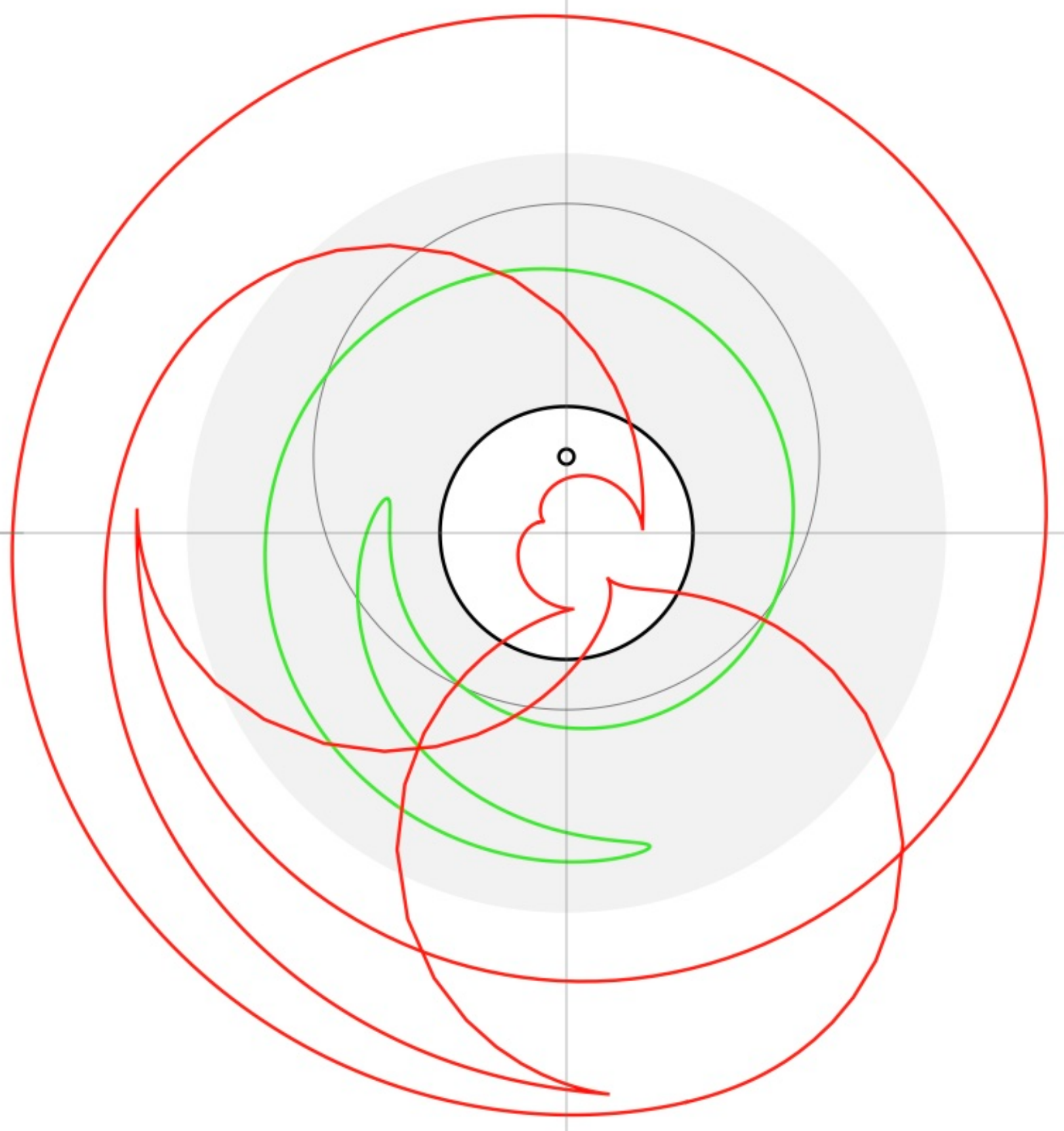}
\caption{Same billiard with a different choice of the initial point $O$ located farther from the center of the disc.}
\label{circles1}
\end{figure}

For example, consider a constant weak magnetic field in an oval. The billiard trajectories are arcs of a circle of radius $R$, and the weakness of the field means that the minimal curvature of the oval is greater than $1/R$ (hence the trajectories cannot touch the oval from inside).
The caustic by reflection is the envelope of a 1-parameter family of circles of radius $R$, and it has two components. See Figure \ref{circles} and \ref{circles1} for the first and second caustics by reflection in a circle. 

\paragraph{Acknowledgment.} I am grateful to Gil Bor for numerous interesting discussions and help with computer experiments. 
This research was supported by NSF grant DMS-2404535 and by Simons grant TSM-00007747.

\end{document}